\documentclass[english]{amsart}
\usepackage[T1]{fontenc}
\usepackage[latin9]{inputenc}
\usepackage{units}
\usepackage{amsmath}
\usepackage{amssymb}
\usepackage{nicefrac}
\usepackage{units}
\usepackage{babel}
\usepackage{geometry}
\usepackage{tkz-graph}
\usepackage{bbm}
\usepackage{a4wide}
\usepackage{graphicx}
\usepackage{color}

\makeatletter

\newcommand{\mathcircumflex}[0]{\mbox{\^{}}}

\makeatother

\begin{document}

\author[Edinah K. Gnang]{Edinah K. Gnang\,$^\ast$}
\address{Edinah K. Gnang, Computer Science Department, Rutgers University (New Brunswick), Piscataway, NJ, USA.}
\email{gnang@cs.rutgers.edu}

\author[Patrick Devlin]{Patrick Devlin}
\address{Patrick Devlin, Mathematics Department, Rutgers University (New Brunswick), Piscataway, NJ, USA.}
\email{prd41@math.rutgers.edu}

\thanks{$^{\ast}$ Supported by the NSF-DGE-0549115 grant}

\date{}

\title{ Some integer formula-encodings and related algorithms }

\date{}

\maketitle

\begin{abstract}
We investigate the special class of formulas made up of arbitrary but finite combinations of 
addition, multiplication, and exponentiation gates. The inputs to these formulas are restricted 
to the integral unit $1$. In connection with such formulas, we describe two essentially distinct families 
of canonical formula-encodings for integers, respectively deduced from the decimal encoding and the 
fundamental theorem of arithmetic. Our main contribution is the detailed description of two algorithms 
which efficiently determine the canonical formula-encodings associated with relatively large sets of 
consecutive integers.
\end{abstract}

\section{Introduction}
It is a well known fact that the binary encoding is on average optimal for representing integers.
However if we think of a binary string as a computer program, it follows that such a program 
implicitly describes a circuit which evaluates to the corresponding integer. When literaly interpreted, 
the binary representation describes a sum of powers of two with the powers determined by the 
location of the bits. The recursive encoding which explicitly describes the 
circuit representation associated with the literal interpretation of the decimal strings was 
pioneered by Goodstein \cite{Gd}. 
In the current discussion we depart slightly from conventional arithmetic circuit models \cite{K,MF} in the 
fact that we consider circuits or more specifically formulas which combine fan-in two exponentiation, 
multiplication, and addition gates with input restricted to the integral unit $1$.
While at first it might seem unnatural to allow exponentiation gates, we argue that exponentiation 
gates are implicit in the decimal encoding. Furthermore, exponentiation 
gates are critical for obtaining a small circuit which evaluate to the integer specified by the input binary strings. 
Throughout the discussion the arithmetic formulas will be described with symbolic expressions and
for convenience we associate with the symbol $x$ the recurring formula $\left(1+1\right)$.
Our main contribution is an asymptotically optimal algorithm for finding Goodstein formula 
encodings for relatively large subset of consecutive integers. Finally we describe an alternative 
canonical formula-encoding conjecturely smaller on average when compared with the Goodstein 
formula-encoding. We also provide an efficient algorithm for computing the latter formula-encodings
for relatively large subsets of consecutive integers.

\section{The Set of Formula-Encodings of Positive Integers}
Let $\mathcal{E}$ denote the set of symbolic expressions
which result from finite combinations of addition, multiplication 
and exponentiations where the {\it only} input is $1$.
For instance  (abbreviating, for convenience, $x=1+1$)
\begin{equation}
x\mathcircumflex\left(1+\left(x\cdot x\cdot x\right)+x\mathcircumflex\left(x\cdot x\right)\right)+1\mathcircumflex\left(1\mathcircumflex1\right)+1\mathcircumflex\left(x\mathcircumflex1\right)+1\mathcircumflex\left(x\mathcircumflex\left(x\mathcircumflex1\right)\right)+1\in\mathcal{E}.
\end{equation}
Elements of the set $\mathcal{E}$ for our
purposes will be encoded as strings from the alphabet $\mathfrak{A}$
\begin{equation}
\mathfrak{A}\;:=\left\{ 1,\: \:+,\:\cdot,\:\mathcircumflex\right\} .
\end{equation}

For the reader's convenience we shall adopt the infix notation thereby making
use of the parenthesis characters '$($' and '$)$'. However we point out that
the parenthesis characters '$($' and '$)$' can be omitted from the alphabet $\mathfrak{A}$
since either the postfix or prefix notations avoid their use entirely.\\

Of course every such expression evaluates to a positive integer, and
we are interested in the shortest possible expression of representing
any given positive integer, or at least (for large integers) 
as close as possible.
The evaluation function is defined recursively as
\begin{equation}
E(a+b)=E(a)+E(b) \quad, \quad
E(a\cdot b)=E(a)\cdot E(b) \quad, \quad
E(a \mathcircumflex b)=E(a)^{E(b)}.
\end{equation}
One can introduce ``axioms'' that transform one tree
to another without changing its value, but these are left to the reader.
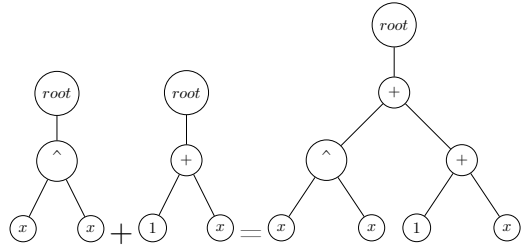
\begin{figure} 
\center
\scalebox{0.6}{
\begin{tikzpicture}[level/.style={sibling distance=30mm/#1}]
\node [circle,draw] (z){$root$}
  child [grow = down] {node [circle,draw] (b) {\Huge $\mathcircumflex$}
  	child {node [circle,draw] (c) {$x$}}
  	child {node [circle,draw] (d) {$x$}}
  	};	
\end{tikzpicture}

{\Huge $+$}

\begin{tikzpicture}[level/.style={sibling distance=30mm/#1}]
\node [circle,draw] (z){$root$}
  child [grow = down] {node [circle,draw] (b) {$+$}
  	child {node [circle,draw] (c) {$1$}}
  	child {node [circle,draw] (d) {$x$}}
  	};	
\end{tikzpicture}
{\Huge =}
\begin{tikzpicture}[level/.style={sibling distance=60mm/#1}]
\node [circle,draw] (y){$root$}
child{ node [circle,draw] (z){$+$}
    child {node [circle,draw] (b) {\Huge $\mathcircumflex$}
    	    child {node [circle,draw] (c) {$x$}}
  	    child {node [circle,draw] (d) {$x$}}
    }
    child {node [circle,draw] (b) {$+$}
            child {node [circle,draw] (c) {$1$}}
  	    child {node [circle,draw] (d) {$x$}}
        }
    };
\end{tikzpicture}
}
\caption{Illustration of a formula operation.}
\end{figure}

\section{Canonical forms}
We shall crucially require for our purposes the notion of 
\emph{canonical form expressions}  or 
\emph{canonical form formulas}.
 Canonical form expressions or formulas are elements of 
$\mathcal{E}$
which we think of as unambiguous representatives of 
the corresponding integer.
We will discuss here two important canonical 
forms.
We point out however that our choices of canonical forms are
bound to be somewhat arbitrary and incidentally alternative representative
choices could be made. 

\subsection{The First Canonical Form}
An expression $f\in\mathcal{E}$ is in the First Canonical Form (FCF) if $f$ corresponds to a finite
sum of the form  (recall that $x$ is short for $1+1$)
\begin{equation}
f=\sum_{k}\left(x\mathcircumflex f_{k}\right)\mbox{ or }f=1+\sum_{k}\left(x\mathcircumflex f_{k}\right)
\end{equation}
such that the expressions $f_{k}$ are distinct for distinct values
of the index $k$ and each ones of the expressions 
$f_{k}\in\mathcal{E} $
being themselves in the FCF.\\
 \\
\textbf{\emph{Proposition 1:}} An arbitrary 
$f\in\mathcal{E}$ is either in the FCF or can be transformed into an expression in the
FCF via a finite sequence of transformations which preserve evaluation
value.\\
 \\
\textbf{\emph{Proposition 2:}} Every expression 
$f\in\mathcal{E}$,
is a member of finite set of non trivial equivalent expressions (i.e. expressions not including any subexpressions of the form $g\mathcircumflex 1$, $1\mathcircumflex g$ or $g \cdot 1$ for some arbitrary expression $g$).\\
\textbf{\emph{Proof : }} A constructive proof of proposition 1 and 2 readily follows from the
quotient remainder theorem.

\subsection{The Second Canonical Form}
An expression $f\in\mathcal{E}$
is in the Second Canonical Form (SCF) if $f$ corresponds to a finite
product of the form 
\begin{equation}
f=\prod_{k}\left[\left(1+f_{k}\right)\mathcircumflex g_{k}\right]\mbox{ or }f=\left(x\mathcircumflex g\right)\cdot\prod_{k}\left[\left(1+f_{k}\right)\mathcircumflex g_{k}\right]
\end{equation}
where for distinct values of the index $k$, the formula associated with $\left(1+f_{k}\right)$
encodes distinct primes greater than $2$.
Furthermore, $f_{k},\: g_{k},\: g\in\mathcal{E}$ are themselves expressions 
in the SCF.\\
\\
\textbf{\emph{Proposition 3:}} An arbitrary 
$f\in\mathcal{E}$
is either in the SCF or can be transformed into an expression in
the SCF via a finite sequence of transformations which preserve the 
evaluation value.\\
\textbf{\emph{Proof : }} A constructive proof of proposition 3 immediately follows from the
fundamental theorem of arithmetic.\\
A considerable advantage of the SCF as a default encoding for integers
is the fact that the encoding considerably simplifies the computational 
complexity analysis of formulas arithmetic. In particular the complexity 
analysis of formulas arithmetic (mulitiplication and exponentiation) reduces 
to the analysis of the formula addition operations. Furthermore it has been 
empirically observed that the lengths of expressions describing formulas in 
the SCF have smaller expected length than their FCF counterpart. 
We further remark that the SCF is implicit in the discussion of integer 
prime tower encodings \cite{GT}.

\section{Computing FCF integer encodings}
We describe here an asymptotically optimal algorithm for determining symbolic expressions which 
describe FCF integer formula-encoding for relatively large set of consecutive integers. 
The algorithm is based on the observation that given the FCF encoding of the first $n$
positive integers one easily deduces from them the FCF encoding for the
next $2^n-n$ positive integers. We pointed out earlier that the FCF encoding 
describes formulas corresponding to the Goodstein base $2$ recursive (or hereditary) integer 
encoding \cite{Gd}, however it is clear that an attempt to uncover the FCF encoding of 
integers by simply iterating through consecutive integers and recursively expressing in 
binary form the powers of $2$, would yield a very inefficient algorithm.
Incidentally our proposed algorithm for determining FCF integer formula-encodings 
amounts to a set recurrence. The initial sets for the recurrence is specified by 
\[
\mathbb{N}_{0} := \left\{ 1 \right\}
\]
and the set recursion is defined by 
\[
\mathbb{N}_{k+1} = \bigcup_{S \in \left\{ \left\{ 1 \right\} \: \cup \: x^{\mathbb{N}_k} \right\}} \left\{ \sum_{s \in S} s \right\}. 
\]
where for an arbitrary symbolic expression $f$ and a set of symbolic expressions $L$, the set $f^L$ is to be 
interpreted as: 
\[
\left\{f^l\right\}_{\left\{l\in L\right\}}. 
\]
Three iterations of the set recurrence yield 
\[
\mbox{FCF}_3 :=
\]
\[
 \left\{ 1, \, x, \, (x+1), \, x^x,  \, (x^x+1), \cdots , \, ( x^{x} + x^{\left(x^{x}\right)} + x + x^{\left(x^{x} + x + 1\right)} + x^{\left(x^{x} + x\right)} + x^{\left(x^{x} + 1\right)} + x^{\left(x + 1\right)} + 1) \right\}.
\]

It follows from the definition of the set recurrence, that the proposed algorithm requires 
$O\left(\log_{\star}(n)\right)$ iterations to produce FCF formulas for all positive integers less 
than $n$ and the algorithm requires optimally $O\left(n \right)$ symbolic expression manipulations.

\section{Zeta recursion.}
We recall here an elementary recursion called the \emph{Zeta recursion} 
emphasizing its close resemblance with the Zeta summation formula. 
Let us briefly recall here the \emph{Zeta recursion} first introduced in \cite{GT} as a combinatorial 
construction for sifting primes.
\begin{equation}
\mathbb{P}_{0}\,:=\left\{ x\right\} ,\:\check{\mathbb{N}}_{0}\,:=\mathbb{P}_{0}\cup\left\{ 1\right\} 
\end{equation}
we consider the set recurrence relation defined by
\begin{equation}
\mathbb{N}_{k+1}=\prod_{p\in\mathbb{P}_{i}}\left\{ \left\{ 1\right\} \cup p^{\check{\mathbb{N}}_{k}}\right\} ,
\end{equation}
where 
\[
p^{\check{\mathbb{N}}_{k}} :=\left\{p^n \mbox{ such that } n \in \check{\mathbb{N}}_{k} \right\}, 
\]
and for sets of symbolic expressions $\left\{S_i\right\}_{0 \le i<m}$
\[
\prod_{0 \le i<m}S_i := \left\{ \prod_{0 \le k < m} s_k \right\}_{s_i \in S_i}.
\]
Finally $\check{\mathbb{N}}_{k+1}$ is deduced from $\mathbb{N}_{k+1}$ by adjunction of missing primes suggested 
by identification of gaps of size two between consecutive elements of $\mathbb{N}_{k+1}$, hence 
\begin{equation}
\mathbb{P}_{k+1}=\mathbb{P}_{k}\cup\left(\check{\mathbb{N}}_{k+1}\backslash\mathbb{N}_{k+1}\right),
\end{equation}
so that for all $k\ge0$, we have $\mathbb{P}_{k}\subsetneqq\check{\mathbb{N}}_{k}$
and 
\begin{equation}
\begin{cases}
\begin{array}{c}
\check{\mathbb{N}}_{k}\subsetneqq\check{\mathbb{N}}_{k+1}\\
\mathbb{P}_{k}\subsetneqq\mathbb{P}_{k+1}
\end{array}\end{cases}.
\end{equation}
Furthermore we can use the \emph{Zeta recursion} to iteratively construct larger 
and larger subsets of rational numbers deducing the set $\mathbb{Q}_{k}$ 
from the previously obtained sets $\check{\mathbb{N}}_{k}$ and $ $$\mathbb{P}_{k}$
as follows 
\begin{equation}
\mathbb{Q}_{k}=\prod_{p\in\mathbb{P}_{k}}\left\{ \left(\frac{1}{p}\right)^{\check{\mathbb{N}}_{k}}\cup\left\{ 1\right\} \cup p^{\check{\mathbb{N}}_{k}}\right\} ,
\end{equation}
where
\begin{equation}
\mathbb{Q}_{k}\subsetneqq\mathbb{Q}_{k+1}.
\end{equation}
Which yields an alternative combinatorial proof of Cantor's result establishing that 
the rational numbers are countable. 

\subsection{Improved Zeta recursion}
Some slight modifications to the \emph{Zeta recursion} has the benefit of improving the 
computational performance of the recurrence computation. 
\begin{equation}
\mathbb{P}_{0}\,:=\left\{ x\right\} ,\:\check{\mathbb{N}}_{0}\,:=\mathbb{P}_{0}\cup\left\{ 1\right\} 
\end{equation}
for an arbitrary $q\in\mathbb{P}_{k}$ we have 
\begin{equation}
\mathbb{N}_{q,k+1}=\bigcup_{\begin{array}{c}
n\in\check{\mathbb{N}}_{k}\\
q^{n}<2^{k+2}
\end{array}}\left\{ \left[2^{k+1},\:2^{k+2}\right]\cap\left(q^{n}\times\prod_{\begin{array}{c}
p\in\mathbb{P}_{k}\\
p<q
\end{array}}\left\{ \left\{ 1\right\} \cup p^{\check{\mathbb{N}}_{k}}\right\} \right)\right\} 
\end{equation}
from which we have that 
\begin{equation}
\mathbb{N}_{k+1}=\check{\mathbb{N}}_{k}\cup\left(\bigcup_{q\in\mathbb{P}_{k}}\mathbb{N}_{q,k+1}\right).
\end{equation}
The completion of the set $\mathbb{N}_{k+1}$ to $\check{\mathbb{N}}_{k+1}$
is still determined by sorting the element in the set $\bigcup_{q\in\mathbb{P}_{k}}\mathbb{N}_{q,k+1}$
and adjoining missing primes located by identifying gaps of size two between consecutive elements
of $\mathbb{N}_{k+1}$.\\
The improved Zeta recursion is not a particularly efficient algorithm for the sole purpose of siefting primes 
because it implicitly requires us to store rather large list of integers. The algorithm is however particularly 
well suited to the task of determining symbolic expressions describing SCF formula-encoding for a relatively 
large set of consecutive integers, with initial sets for the iteration being 
\[
\mathbb{P}_{0}\,:=\left\{ x\right\} ,\:\check{\mathbb{N}}_{0}\,:=\mathbb{P}_{0}\cup\left\{ 1\right\}. 
\]
For instance seven iterations of the improved Zeta recursion yield
\[
\mbox{SCF}_7  =  \left\{1,\, x, \, (x+1),\, x^x, \, (x^x+1),\, \left(x + 1\right) x ,\, \left(\left(x + 1\right) x + 1\right),\, \cdots,\, \left(x + 1\right) \left(x^{\left(x^{x}\right)} + 1\right) \left(x^{x} + 1\right) \right\} 
\]
It follows from the definition of the Zeta recursion that the proposed algorithm requires $O\left(\log n\right)$
iterations to determine SCF formulas for all positive integers less than $n$ and the algorithm requires 
$O\left(\frac{n^2}{\log n} \right)$ symbolic expression manipulations.

\section{Acknowledgments}
We are especially grateful to Professors Doron Zeilberger, Henry Cohn, and Mario Szegedy for insightful comments while
preparing this manuscript.

\end{document}